\documentclass[english,a4paper,twoside,11pt,notitlepage]{amsart}

\usepackage[T1]{fontenc}
\usepackage[latin1]{inputenc}
\usepackage{babel,amsfonts,amsmath,amsthm,float,amscd,amssymb,latexsym}
\usepackage[dvips]{epsfig}
\usepackage{fullpage}
\usepackage{xypic}
\usepackage{url}

\theoremstyle{plain}
\newtheorem{theorem}{Theorem}
\newtheorem{lemma}[theorem]{Lemma}

\newtheorem{conjecture}[theorem]{Conjecture}
\newtheorem{corollary}[theorem]{Corollary}

\theoremstyle{definition}
\newtheorem{definition}[theorem]{Definition}

\newtheorem{problem}[theorem]{Problem}
\newtheorem{proposition}[theorem]{Proposition}

\newcommand{\di}{\mathrm{d}}

\newcommand{\gt}{\tilde{g}}
\newcommand{\an}{\frac{n-2}{4(n-1)}}

\newcommand{\Nn}{\frac{2n}{n-2}}

\newcommand{\e}{\varepsilon}

\title{The Yamabe problem with singularities}

\author{Farid Madani}
\address{Institut Math\'ematiques de Jussieu, Universit\'e Pierre et Marie Curie\\
\'Equipe d'Analyse Complexe et G\'eom\'etrie \\
175, rue Chevaleret\\
75013 Paris, France.}

\date{}
\begin{document}

\maketitle

\begin{abstract}
Let $(M,g)$ be a compact Riemannian manifold of dimension $n\geq 3$. 
Under some assumptions, we prove that there exists a positive function  $\varphi$ solution of the following Yamabe type equation
\begin{equation*}
\Delta \varphi+ h\varphi= \tilde h \varphi^{\frac{n+2}{n-2}}
\end{equation*}
where $h\in L^p(M)$, $p>n/2$  and $\tilde h\in \mathbb R$. We give the regularity of $\varphi$ with respect to the value of $p$. 
Finally, we consider the results in geometry when $g$ is a singular Riemannian metric and $h=\frac{n-2}{4(n-1)}R_g$,  where $R_g$ is the scalar curvature  of $g$.
\end{abstract}


\section{Introduction}
Let  $(M,g)$ be a smooth compact Riemannian manifold of dimension $n\geq 3$. Denote by $R_g$ the scalar curvature of $g$. The Yamabe problem is the following: 
\begin{problem}\label{yamabepro}
Does there exists a constant scalar curvature metric conformal to $g$? 
\end{problem}

If  $\tilde g=\varphi^{4/(n-2)}g$ is a conformal metric to $g$ with $\varphi$ a smooth positive  function, then the scalar curvatures  $R_g$ and $R_{\gt}$ are related by the following equation:

\begin{equation}\label{yamabe}
\frac{4(n-1)}{n-2}\Delta_g \varphi+ R_g\varphi=  R_{\gt} \varphi^{N-1}
\end{equation}
where $N=\Nn$ and $\Delta_g$ is the geometric Laplacian of the metric $g$ with nonnegative eigenvalues.\\ 
To solve the Yamabe problem, it is equivalent to find a function $\varphi$ solution of equation above where $R_{\gt}$ is constant. Equation \eqref{yamabe} is called Yamabe equation. Yamabe \cite{Yam} stated the following functional, defined for any $\psi \in H_1(M)-\{0\}$ by
\begin{equation}\label{fonctionnel}
I_g(\psi)=\frac{E(\psi)}{\|\psi\|_N^2}=\frac{\displaystyle\int_M |\nabla\psi|^2+
\frac{n-2}{4(n-1)}R_g\psi^2\di v}{\|\psi\|_N^2}
\end{equation}
and he considered the infimum of  $I_g$ defined as follow
$$\mu(g)=\inf_{\psi\in H_1(M)-\{0\}}I_g(\psi)$$
He solved the case when $\mu(g)$ is nonpositive.  Aubin \cite{Aub} showed that it was sufficient to solve the following conjecture:
\begin{conjecture}[Aubin \cite{Aub}]\label{Aubincon}
If $(M,g)$ is not conformal to $(S_n,g_{can})$ then
 \begin{equation}\label{cau}
 \mu(M,g)<\mu(S_n,g_{can})
\end{equation}
where $\mu(M,g)=\inf\{I_g(\psi),\; \psi\in H_1(M)-\{0\}\}$ 
\end{conjecture}
It is known that $\mu(S_n,g_{can})=K^{-2}(n,2)=\frac{1}{4}n(n-2)\omega_n^{2/n}$, where $\omega_n$ is the volume of the unit sphere $S_n$ and $K(n,2) $ is defined in  theorem \ref{Best c}.\\
In the following, we write  $\mu(g)$ instead of $\mu(M,g)$.\\
Aubin proved that the conjecture is valid for all smooth compact non conformally flat Riemannian manifolds of dimension  $n\geq 6$ and conformally flat manifolds with finite  non trivial fundamental group. The case of conformally flat manifolds and the dimensions  3,4 and 5 were solved by   Schoen \cite{Schoen} using positive mass theorem. Hence the conjecture above holds. By works of Yamabe \cite{Yam}, Aubin \cite{Aub} and Schoen \cite{Schoen}, the Yamabe problem is completely solved, when the manifold is compact and smooth. \\
The purpose of this paper is to study the following equation
\begin{equation}
 \Delta_g\psi+h\psi=\tilde h\psi^{\frac{n+2}{n-2}}
\end{equation}
where $h\in L^p(M)$, and $\tilde h\in \mathbb R$. We call this kind of equation "Yamabe type equation".  We will give a  special consideration for the case  $h=\frac{n-2}{4(n-1)}R_g$.

\section{Regularity theorems for Yamabe type equations}

\begin{theorem}\label{reg}
Let $\Omega$ be an open subset of $\mathbb{R}^n$ and $L$ an uniformely elliptic  linear operator of the second degree, defined by
\begin{equation}\label{elli}
L(u)=\sum_{i,j}a_{ij}\partial_{ij}u+ \sum_i b_i\partial_iu+hu
\end{equation}
where the coefficients $a_{ij},\;b_i\mbox{ and }h $ are real valued bounded functions of class $C^k$ with $k\in \mathbb N$. Let $u$ be a weak solution  of the equation $Lu=f$. 
\begin{enumerate}
\item[(i)] If $f\in C^{k,\alpha}(\Omega)$ then $u\in C^{k+2,\alpha}(\Omega)$
\item[(ii)]If $f\in H_k^p(\Omega)$ then $u\in H_{k+2}^p(\Omega)$
\end{enumerate}
\end{theorem}
This theorem is the standard regularity theorem, we can find a proof in the book of Gilbarg and Trudinger~\cite{GT}.\\
The following two theorems allow us to find the best regularity for the solution of Yamabe type equations. Using the theorem \ref{gef}, Trudinger \cite{Trud} Showed that the weak solutions of the Yamabe equation \eqref{yamabe} are smooth. Yamabe \cite{Yam} had already used implicitely this theorem.

\begin{theorem}\label{gef}
 On a $n-$dimensional compact Riemannian manifold $(M,g)$, if $u\geq 0$ is a non trivial weak solution in $H_1(M)$ of  equation $\Delta_g u+hu=0$, with $h\in L^p(M)$ and $p>n/2$, then $u\in C^{1-[n/p],\beta}(M)$  and positive.\\
$[n/p]$ is the integer part of $n/p$,  $\beta\in(0,1)$.
\end{theorem}

Notice that if $u$ satisfies the assumptions of this  theorem, then $\Delta u\in L^p(M)$. Regularity theorem \ref{reg}  implies that $u\in H^p_2(M)$ and using Sobolev embedding, we find $u\in C^{1-[n/p],\beta}(M)$\\
Theorem \ref{gef} permits  to proof the following theorem:
\begin{theorem}\label{regYS}
Let $(M,g)$ be $n-$dimensional compact smooth Riemannian manifold. $p$ and $\tilde h$ are two reel numbers, with $p>n/2$. If $\varphi\in H_1(M)$ is a non trivial, nonnegative weak solution of  
\begin{equation}\label{etyam}
 \Delta_g\psi+h\psi=\tilde h\psi^{\frac{n+2}{n-2}}
\end{equation}
 
then $\varphi\in H_2^p(M)\subset C^{1-[n/p],\beta}(M)$ and  $\varphi$ is positive.
\end{theorem}

\begin{proof}
It is sufficient to show that there  exists $\e>0$ such that $\varphi\in L^{(\e+2n)/(n-2)}(M)$. Indeed, if $\varphi$ satisfies the assumptions of theorem and  belongs to $L^{(\e+2n)/(n-2)}(M)$, then it is a solution of
$$\Delta_gu+(h-\tilde h\varphi^{\frac{4}{n-2}})u=0$$
with $h-\tilde h\varphi^{\frac{4}{n-2}}\in L^r(M)$ and $r=\min(p,\frac{2n+\e}{4})>n/2$. Using theorem \ref{gef}, we deduce  that $\varphi$ is positive and continuous. Theorem \ref{reg} and Sobolev embedding imply that $\varphi\in H^p_2(M)$ with $p>n/2$.\\ 
Let $l$ be a positive reel number and  $H$, $F$ are two continuous functions  in $\mathbb R_+$ defined by:
\begin{align*}
 H(t) &=\begin{cases} 
t^\gamma &\text{ if } 0\leq t\leq l\\
l^{q-1}(ql^{q-1}t-(q-1)l^q) &\text{ if }  t>l
      \end{cases}\\
F(t) &=\begin{cases}
t^q &\text{ if } 0\leq t\leq l\\
 ql^{q-1}t-(q-1)l^q \qquad &\text{ if }  t>l
     \end{cases}
\end{align*}

\begin{equation}
\text{where }\gamma=2q-1,\text{ and } 1<q<\frac{n(p-1)}{p(n-2)}  
\end{equation}
 $\varphi$ is positive, belongs to $ H_1(M)$. $H\circ\varphi$ and $F\circ \varphi$ belong also to $H_1(M)$. Notice that for any $t\in\mathbb R_+ -\{l\}$
\begin{equation}\label{forl}
 qH(t)=F(t)F'(t),\; (F'(t))^2\leq qH'(t)\text{ and }F^2(t)\geq tH(t)
\end{equation}
If $\varphi $ is a weak solution of equation \eqref{etyam}, then
\begin{equation}\label{tyfai}
\forall\psi\in H_1(M)\quad \int_M\nabla\varphi\cdot\nabla\psi \di v+\int_Mh\varphi\psi\di v=\tilde h\int_M\varphi^{N-1}\psi \di v
\end{equation}
where $N=2n/(n-2)$.\\
Let us choose $\psi=\eta^2H\circ\varphi$, where $\eta$ is $C^1-$function with support in the ball $B_P(2\delta)$ and radius $2\delta$ sufficiently small, such that $\eta=1$ on $B_P(\delta)$. If we substitute in \eqref{tyfai}, we obtain
\begin{equation}\label{refref}
 \int_M \eta^2H'\circ\varphi|\nabla\varphi|^2\di v+2\int_M\eta H\circ\varphi\nabla\varphi\cdot\nabla\eta\di v=\tilde h\int_M\varphi^{N-1}\eta^2 H\circ\varphi\di v-\int_Mh\varphi\eta^2H\circ\varphi\di v
\end{equation}
Let $f=F\circ\varphi$ be a function. We estimate the forth integrals above, using  function $f$ and relations \eqref{forl}. We have $\nabla f=F'\circ\varphi\nabla\varphi$, the second relation in \eqref{forl} implies
$$|\nabla f|^2=(F'\circ\varphi)^2|\nabla\varphi|^2\leq qH'\circ\varphi|\nabla\varphi|^2$$ 
We deduce that the first integral of equality \eqref{refref} is bounded from below. 
$$\frac{1}{q}\|\eta \nabla f\|_2^2\leq\int_M \eta^2H'\circ\varphi|\nabla\varphi|^2\di v$$
The first relation of \eqref{forl} and Cauchy--Schwarz inequality imply that the second integral of \eqref{refref} is bounded from below by:
$$2\int_M\eta H\circ\varphi\nabla\varphi\cdot\nabla\eta\di v=\frac{2}{q}\int_M\eta f\nabla f\nabla\eta\di v
\geq \frac{-2}{q}\|f\nabla\eta\|_2\|\eta\nabla f\|_2$$
By the last relation in \eqref{forl}, we have $\varphi H\circ\varphi\leq f^2$. The two integrals in the right side in \eqref{refref} are bounded  by:
$$\biggl|\tilde h\int_M\varphi^{N-1}\eta^2 H\circ\varphi\di v-\int_Mh\varphi\eta^2H\circ\varphi\di v\biggr|\leq |\tilde h|\|\varphi\|^{4/(n-2)}_{N,2\delta}\|\eta f\|^2_N+\|h\|_p\|\eta f\|_{2p/(p-1)}^2$$  
where $\|\varphi\|^N_{N,r}=\int_{B_P(r)}\varphi^N\di v$. If we take together these estimates, equality \eqref{refref} becomes:
\begin{equation}\label{refre}
 \|\eta \nabla f\|_2^2-2\|f\nabla\eta\|_2\|\eta\nabla f\|_2\leq q(|\tilde h|\|\varphi\|^{4/(n-2)}_{N,2\delta}\|\eta f\|^2_N+\|h\|_p\|\eta f\|_{2p/(p-1)}^2)
\end{equation}
Notice that for all nonnegative reel numbers $a,\;b,\;c\text{ and }d$, if $a^2-2ab\leq c^2+d^2$ then $a\leq c+d+2b$. Using this remark, inequality \eqref{refre} becomes:
\begin{equation}\label{inega Est}
 \|\eta \nabla f\|_2\leq \sqrt{q|\tilde h|}\|\varphi\|^{2/(n-2)}_{N,2\delta}\|\eta f\|_N+\sqrt{q\|h\|_p}\|\eta f\|_{2p/(p-1)}+2 \|f\nabla\eta\|_2
\end{equation}
By Sobolev embedding, we know that there exists a positive constant $c$, which depends only on  $n$, such that 
$$\|\eta f\|_N\leq c(\|\eta \nabla f\|_2+\|f\nabla\eta\|_2+\|\eta f\|_2)$$
The choice of $q$ ($q< N$) and inequality \eqref{inega Est} permit to write  
\begin{equation}
  (1-c\sqrt{N|\tilde h|}\|\varphi\|^{2/(n-2)}_{N,2\delta})\|\eta f\|_N\leq c\bigl(\sqrt{N\|h\|_p}\|\eta f\|_{2p/(p-1)}+3 \|f\nabla\eta\|_2+\|\eta f\|_2\bigr)
\end{equation}
We choose $\delta$ sufficiently small such that 
$$\|\varphi\|^{2/(n-2)}_{N,2\delta}\leq 1/(2c\sqrt{N|\tilde h|})$$ 
when  $l$ goes to $+\infty$, we deduce that there exists a postive constant $C$, which depends on $n,\;\delta, \|\eta\|_\infty,\;\|\nabla\eta\|_\infty,\;\|h\|_p$ and $|\tilde h|$  such that
$$\|\varphi^q\|_{N,2\delta}\leq C(\|\varphi^q\|_2+\|\varphi^q\|_{2p/(p-1)})$$
$\frac{2p}{p-1}q<N$ and $\varphi$ is bounded in $L^N$, hence 
$$\|\varphi\|_{qN,2\delta}\leq C$$ 
If $(\eta_i)_{i\in I}$ is a partition of unity subordinate to the  covering  $\{B_{P_i}(\delta)\}_{i\in J}$ on $M$ $$\|\varphi\|^{qN}_{qN}=\sum_{i\in I}\|\eta_i\varphi\|^{qN}_{qN,\delta_i}\leq C$$ 
Hence $\varphi\in L^{qN}$ with $qN>N$. The remark in the begining of the proof implies the theorem.
\end{proof}

\begin{proposition}\label{delta fu}
Let $(M,g)$ be $n-$dimensional compact smooth Riemannian manifold. $L:=\Delta_g+h$ is a linear operator, with $h\in L^p(M)$ and $p>n/2$. If the smallest eigenvalue $\lambda$ of $L$  is positive then
\begin{itemize}
 \item[i.] $L$ est coercive, in other words there exists $c>0$ such that
$$\forall\psi\in H_1(M)\quad (L\psi,\psi)_{L^2}\geq c(\|\nabla\psi\|^2_2+\|\psi\|^2_2)$$
\item[ii.] The opertor $L: H_2^{p}(M)\longrightarrow L^p(M)$ is invertible. 
\end{itemize}

\end{proposition}

\begin{proof} $L$ admits a smallest eigenvalue because if $\lambda$ is an eigenvalue associated to the eigenfunction $\psi$ then there exists $C>0$ such that
$$\lambda\|\psi\|_2^2=(L\psi,\psi)_{L^2}=\|\nabla\psi\|_2^2+\int_Mh\psi^2\di v\geq -\|h\|_p\|\psi\|_{2p/(p-1)}^2\geq -C\|h\|_p\|\psi\|_2^2$$

Hence $\lambda\geq -C\|h\|_p$. If $\lambda$ is the smallest eigenvalue of $L$ then
$$\lambda=\inf_{\varphi\in H_1(M)-\{0\}}\frac{E(\varphi)}{\|\varphi\|_2^2}$$
where $$E(\varphi)=(L\varphi,\varphi)_{L^2}=\int_M|\nabla\varphi|^2+h\varphi^2\di v$$ 
So, for any $\varphi\in H_1(M)$
\begin{equation}\label{vp lam}
E(\varphi)\geq \lambda\|\varphi\|_2^2
\end{equation}

Suppose that $L$ is non coercive, then there exists a sequence $(\psi_i)_{i\in\mathbb N}$ in $H_1(M)$, which satisfies
$$E(\psi_i)<\frac{1}{i}(\|\nabla\psi_i\|^2_2+vol(M)^{2/n})\mbox{ and }\|\psi_i\|_N=1$$
It implies
$$(1-\frac{1}{i})E(\psi_i)<\frac{vol(M)^{2/n}}{i}-\frac{1}{i}\int_Mh\psi_i^2\di v$$
because $|\int_Mh\psi_i^2\di v|\leq \|h\|_{n/2}$, $\lim_{i\to +\infty}E(\psi_i)\leq 0$. On other hands $E(\psi_i)\geq \lambda\|\psi_i\|_2^2$ with $\lambda>0$. Which is impossible.\\
If $L\psi=0$, then, using \eqref{vp lam}, $\varphi=0$. So  $L$ is injective. \\
Let $f\in L^p(M)$. Let us prove that the following equation admit a solution $\psi\in H^p_2(M)$
\begin{equation}\label{equation}
 \Delta\varphi+h\varphi=f
\end{equation}
We  minimize the functional $E$ defined in the begining of the proof. Let define $\mu$ as follow
\begin{equation}
\mu=\inf\{E(\varphi)/\varphi\in H_1(M),\; \int_Mf\varphi\di v=1\}
\end{equation}
and $(\psi_i)_{i\in \mathbb N}$ a sequence in $H_1(M)$ which  minimizes $E$, then $$\lim_{i\to+\infty}E(\psi_i)=\mu\text{ and }\int_Mf\psi_i\di v=1$$
Without loss of generalities, we suppose that for any nonnegative integer $i$, $E(\psi_i)\leq\mu+1$. It implies 
$$c(\|\nabla\psi_i\|^2_2+\|\psi_i\|^2_2)\leq E(\psi_i)\leq\mu+1$$
because $L$ is coercive. We conclude that $(\psi_i)_{i\in \mathbb N}$ is bounded in $H_1(M)$. The Kondrakov theorem and Banach theorem imply that there exists a subsequence $(\psi_j)_{j\in \mathbb N}$ such that
\begin{itemize} 
\item[$*$] $\psi_{j}\rightharpoonup\psi$ weakly in $H_1(M)$
\item[$*$] $\psi_{j}\rightarrow\psi$ strongly in $L^s(M)$ for all $1\leq s<N$ 
\item[$*$] $\psi_{j}\rightarrow\psi$ almost everywhere.
\end{itemize}
Then  $(\psi_j)$ converge strongly in  $L^{2p/(p-1)}(M)$ because $2p/(p-1)<N$. So 
$$\int_Mf\psi\di v=1\text{ and } \int_Mh\psi_j^2\di v\rightarrow \int_Mh\psi^2\di v$$ 
The weak convegence in  $H_1(M)$ and the strong convergence $L^2(M)$ imply 
$$\lim_{j\to +\infty}\|\nabla\psi_j\|_2\geq \|\nabla\psi\|_2$$
We conclude that $E(\psi)\leq \mu$, hence $E(\psi)=\mu$. If we write  Euler--Lagrange equation for $\psi$, we find that it is a weak solution in $H_1(M)$ of equation \eqref{equation}. It remains to prove that $\psi\in H^p_2(M)$. Suppose that $\psi\in L^{s_i}(M)$. Then $hu\in L^{\frac{ps_i}{p+s_i}}(M)$, Hence $\Delta u\in L^{\frac{ps_i}{p+s_i}}(M)$. Regularity  theorem \ref{reg} assures that $u\in H^{\frac{ps_i}{p+s_i}}_2(M)$. We know that $H_2^r(M)\subset L^s(M)$ if $r\leq n/2$ with $s=nr/(n-2r)$, and $H^r_2(M)\subset C^{1-[n/r],\beta}(M)$ if $r>n/2$. These inclusions imply the following results 
\begin{equation*}
\begin{cases}
s_0=N\\
 u \in L^{s_{i+1}}(M)\text{ where } s_{i+1}=\frac{nps_i}{np-(p-2n)s_i}& \mbox{ if }s_i\leq \frac{np}{2p-n}\\
u \in H^p_2(M) & \mbox{ if }s_i> \frac{np}{2p-n}
\end{cases}
\end{equation*}
If there exists $i\in \mathbb N$ such that  $s_i>\frac{np}{2p-n}$, which is equivalent to $\frac{ps_i}{p+s_i}>n/2$ then $u\in C^{0,\beta}(M) $, which implies $\Delta u\in L^p(M)$, hence $u \in H^p_2(M)$. If there exists $i\in \mathbb N$ such that  $s_i=\frac{np}{2p-n}$ then $u\in L^\infty(M)$ and we conclude by regularity theorem that $u \in H^p_2(M)$. Suppose that for any $i\in \mathbb N$, $s_i<\frac{np}{2p-n}$, the sequence $(s_i)_{i\in\mathbb{N}}$ is increasing and bounded from above, it converges to $s=0$ which is impossible. 
\end{proof}

\begin{theorem}\label{Best c}
Let $(M,g)$ be a $n-$dimensional smooth compact Riemannian manifold. For all $\e>0$, there exists $A(\varepsilon)>0$ such that 
$$\forall \varphi\in H_1(M)\quad \|\varphi\|_{N}\leq (K(n,2)+\varepsilon)\|\nabla\varphi\|_2+A(\varepsilon)\|\varphi\|_2$$
where $N=\Nn$ and  $K(n,2)=\frac{2}{\sqrt{n(n-2)}}\omega_n^{-1/n}$ 
\end{theorem}
The inequality of this theorem is a particular case of a more general one.  More further details are given in the Aubin's book \cite{Aubin}. 

\section{Existence theorem}\label{section sans sym}

 We consider the following equation :
\begin{equation}\label{AF}
\Delta_g \psi+ h\psi= \tilde h \psi^{\frac{n+2}{n-2}}
\end{equation}
where $\psi\in H_1(M)$, $h\in L^p(M)$ with $p>n/2$ and $\tilde h$  is a  reel number. As mentionned in the introduction, this kind of equation are called Yamabe type equation. In the particular case when $h=\frac{n-2}{4(n-1)}R_g$, equation \eqref{AF} is the Yamabe equation \eqref{yamabe}. To solve this equation, we use the variational method.\\
We define the energy $E$ of $\psi\in H_1(M)$ by: 
\begin{equation}
 E(\psi)=\int_M |\nabla\psi|^2+h\psi^2\di v
\end{equation}
and we consider the fonctional  $I_g$ defined for all $\psi\in H_1(M)-\{0\}$ by
\begin{equation}
I_g(\psi)=\frac{E(\psi)}{\|\psi\|^2_N}\label{ca} 
\end{equation}
We denote
\begin{equation}
 \mu(g)=\inf_{\psi\in H_1(M)-\{0\},\psi\geq 0}I_g(\psi)=\inf_{\|\psi\|_N=1,\psi\geq 0}E(\psi)
\end{equation}

with $N=\Nn$. the main result of this section is 
\begin{theorem}\label{cg}
If $p>n/2$ and $$\mu(g)<K^{-2}(n,2)$$ 
then equation \eqref{AF} admits a positive solution $\varphi\in H^p_2(M)\subset C^{1-[n/p],\beta}(M)$, which minimizes $I_g$ (i.e. $E(\varphi)=\mu(g)=\tilde h$ and $\|\varphi\|_N=1$). where $\beta\in (0,1)$.
\end{theorem}

To proof this theorem, we need the following lemma, proven by Brezis and Lieb\cite{BL}
\begin{lemma}\label{BLL}
Let $(f_i)_{i\in\mathbb N}$ be a sequence of mesurable functions in  $(\Omega,\Sigma,\mu)$. 
If $(f_i)_{i\in\mathbb N}$ is  uniformely bounded in $L^p$ with $0<p<+\infty$ and $f_i\rightarrow f$ almost everywhere, then
$$\lim_{i\to +\infty}[\|f_i\|_p^p-\|f_i-f\|_p^p]=\|f\|_p^p$$ 
\end{lemma}

\begin{proof}[\textbf{Proof of theorem \ref{cg}}]
We  check that $\mu(g)$ is finite. In fact,  using Hölder inequality, we have 
\begin{equation*} 
E(\psi)\geq -\|h\|_{n/2}\|\psi\|^2_N
\end{equation*} 
we deduce that  $\mu(g)\geq -\|h\|_{n/2}>-\infty$.\\
Let $(\varphi_i)_{i\in \mathbb N}$ be a minimizing sequence:
\begin{equation}\label{AA}
E(\varphi_i)=\mu(g)+o(1),\;\|\varphi_i\|_N=1\mbox{ et }\varphi_i\geq 0 \end{equation}
Applying Hölder inequality again for the equation above, we obtain
\begin{gather*}
\|\nabla\varphi_i\|_2^2\leq \|h\|_{n/2}+\mu(g)+o(1)\\
\|\varphi_i\|_2^2\leq (vol(M))^{2/n} 
\end{gather*}
We conclude that $(\varphi_i)_{i\in \mathbb N}$ is bounded in $H_1(M)$. Without loss of generalities, we suppose that there exists $\varphi\in H_1(M)$ such that
\begin{description}
\item[$*$] $\varphi_i\rightharpoonup \varphi$ weakly in $H_1(M)$  
\item[$*$] $\varphi_i\rightarrow \varphi$ strongly in $L^s(M)$ for any $s\in[1,N)$ 
\item[$*$] $\varphi_i\rightarrow \varphi$ almost everywhere
\end{description}
We deduce that
\begin{equation*}
\int_M |h| |\varphi_i-\varphi|^2\di v\leq \|h\|_p\|\varphi_i-\varphi\|_{2p/(p-1)}^2\rightarrow 0
\text{ strongly because }2p/(p-1)<N
\end{equation*}
Let $\psi_i=\varphi_i-\varphi$, then $\psi_i\rightarrow 0$ weakly in $H_1(M)$, strongly in 
$L^q(M)$ for any $q<N$.\\ 
We have $\|\nabla\varphi_i\|_2^2=\|\nabla\psi_i\|_2^2+\|\nabla\varphi\|_2^2+2\int_M\nabla\psi_i\cdot\nabla\varphi \di v$. Hence
\begin{equation*}
E(\varphi_i)=E(\varphi)+\|\nabla\psi_i\|_2^2+o(1)
\end{equation*}
We know that $E(\varphi)\geq \mu(g)\|\varphi\|_N^2$ by definition of  $\mu(g)$, and $E(\varphi_i)=\mu(g)+o(1)$
 by definition of  $(\varphi_i)_{i\in\mathbb N}$. We conclude
\begin{equation}\label{AE}
\mu(g)\|\varphi\|_N^2+\|\nabla\psi_i\|_2^2\leq \mu(g)+o(1)
\end{equation}
Using lemma \ref{BLL} for $(\varphi_i)_{i\in\mathbb N}$, we obtain
\begin{align}\label{AB}
\|\psi_i\|_N^N &+\|\varphi\|_N^N+o(1)=1\\
\|\psi_i\|_N^2 &+\|\varphi\|_N^2+o(1)\geq 1 \label{AC}
\end{align}
Theorem \ref{Best c} gives 
\begin{equation*}
\|\psi_i\|_N^2\leq (K^2(n,2)+\varepsilon)\|\nabla\psi_i\|_2^2+o(1)
\end{equation*}
Inequality \eqref{AC} becomes
\begin{equation*}
(K^2(n,2)+\varepsilon)\|\nabla\psi_i\|_2^2+\|\varphi\|_N^2+o(1)\geq 1
\end{equation*}
Using the last inequality in \eqref{AE}, we obtain
\begin{equation*}
\mu(g)\|\varphi\|_N^2+\|\nabla\psi_i\|_2^2\leq \mu(g)[(K^2(n,2)+\varepsilon)\|\nabla\psi_i\|_2^2+\|\varphi\|_N^2]+o(1)
\end{equation*}
Finally
\begin{equation*}
[1-\mu(g)(K^2(n,2)+\varepsilon)]\|\nabla\psi_i\|_2^2\leq o(1)
\end{equation*}
If $\mu(g)<K^{-2}(n,2)$, we can choose $\varepsilon$ such that the first factor of this inequality becomes positive. 
We deduce that $(\psi_i)_{i\in\mathbb N}$  converges strongly to zero in 
$H_1(M)$, $\varphi_i\rightarrow \varphi$ strongly in $H_1(M)$ and $L^N(M)$. Hence $I_g(\varphi)=\mu(g)$.\\
We have just found a non trivial solution of the following Yamabe type equation
\begin{equation*}
\Delta \psi+ h\psi= \mu(g) \psi^{N-1}
\end{equation*}
which satisfies $\|\varphi\|_N=1$ and $\varphi\geq 0$. Theorem \ref{regYS} implies  $\varphi\in H^p_2(M)\subset C^{1-[n/p],\beta}(M)$ and $\varphi>0$.
\end{proof}

\section{The choice of the metric}\label{hypG}

From now until the end of this paper, $M$ is a compact smooth manifold of dimension $n\geq 3$. Denote by $T^*M$  the cotangent space of $M$.\\

\textbf{Assumption $\boldsymbol{(H)}$:} \emph{$g$ is a metric in the Sobolev space $H_2^p(M,T^*M\otimes T^*M)$ with $p>n$. There exists a point $P_0\in M$ and $\delta>0$ such that $g$ is smooth in the ball $B_{P_0}(\delta)$.}\\

We can suppose that $g$ is  $C^2$ instead of $C^\infty$ in this ball, but it is not an important point.\\ 
Actually our objectif, in this section is to study the Yamabe problem when the metric $g$ admits a finite number of points with singularities and smooth out side these points. The assumption $(H)$ generalizes this  conditions and define the notion of "singularities".\\
By Sobolev embedding, $H_2^p(M,T^*M\otimes T^*M)\subset C^{1,\beta}(M,T^*M\otimes T^*M)$ for some  $\beta\in (0,1)$. Hence the  metrics which satisfy assumption  $(H)$ are $C^{1,\beta}$. The Christoffels belong to $H^p_1\subset C^\beta(M)$. Riemann curvature tensor, Ricci tensor and scalar curvature are in $L^p$. An example of metric which satisfies assumption $(H)$ is $g=(1+d(P_0,\cdot)^{2-\alpha})^m g_0$ where $g_0$  is a smooth metric, $\alpha\in(0,1)$ and  $d(P_0,\cdot)$ is the  distance function. \\
We obtain many results  which are true for metrics in $H^p_2(M,T^*M\otimes T^*M)$, with  $p>n/2$. In the assumption $(H)$, we add the condition that $p>n$ to have a continuous  Christoffels for $g\in H^p_2(M,T^*M\otimes T^*M) $. The assumption $(H)$ is sufficient to prove the Aubin's conjecture \ref{Aubincon} (cf. theorem \ref{conj aub}), and to construct the Green function of the conformal  Laplacian (cf. section \ref{glc}).\\
We consider the following problem:
\begin{problem}\label{yam sing}
 Let $g$  be a metric which satisfies the assumption  $(H)$. Does there exist a  conformal metric  $\tilde g$   for which the scalar curvature  $R_{\gt}$ is constant ?
\end{problem}

It is clear that if the initial metric $g$ is smooth then the problem above is the Yamabe problem  \ref{yamabepro}, which is completely solved. We will prove that the answer to this problem is positive. The following  proposition tell us that the conformal class of the metrics is well defined when the metrics are in $H^p_2$.
\begin{proposition}
Let $g$ be a metric in $H^p_2$ and $\psi\in H^p_2(M)$ a positive function. If $p>n/2$ then the metric $\gt=\psi^{\frac{4}{n-2}}g$ is well defined, and it is in the same space as $g$.
\end{proposition}

\begin{proof}
Using Sobolev embedding, it is easy to check that $H^p_2(M)$ is an algebra for any $p>n/2$. This proposition  is a consequence of this fact.
\end{proof}

In their paper \cite{LP} about the Yamabe problem, Lee and Parker proved that on every compact  Riemanniann manifold $(M,g)$, there exist a normal coordinates system $\{(U_i,x_i)\}_{i\in I}$ and metric   $g'$ conformal to $g$ such that $\det g'=1+O(|x|^m)$ with $m$ as big as we want. Cao \cite{Cao} and Günther \cite{Gun} proved that we can get $\det g'=1$.
 
\begin{definition}
$\gt$ is a  Cao--Günther metric if it is conformal to $g$  and there exist a coordinates system such that $\det \gt=1$.  
\end{definition}

\begin{theorem}[Cao--Günther]\label{caogun}
Let  $M$ be  $C^{a+2,\beta}$ compact manifold  of dimension $n$ with $a\in\mathbb N$, $\beta\in(0,1)$, $g$ be a $C^{a+1,\beta}-$Riemannian metric , and $P$ be a point in $M$. Then there exists a $C^{a+1,\beta'}-$positive function  $\varphi$ with $\beta'\in (0,\beta)$ such that  $\det(\varphi g)=1$ in a normal coordinates system with origin  $P$. 
\end{theorem}
 
Notice that if the metric $g\in H_2^p(M,T^*M\otimes T^*M)$ with $p>n$ then it belongs to  $C^{1,\beta}$. Hence the manifold  $(M,g)$ admits a Cao--Günther metric. It is not really useful to suppose that the metric is smooth in a ball, for the existence of this kind of metrics.

\section{Conformal Laplacian}

\begin{definition}The conformal Laplacian of Riemannian manifold  $(M,g)$ is the operator $L_g$, defined by :  
$$L_g=\Delta_g+\an R_g$$ 
\end{definition}

It is known that the conformal Laplacian, when $g$ is smooth, is conformally invariant. Actually it verifies \eqref{laplinvfai} strongly. We prove that we have this property even when the metric is in $H^p_2(M,T^*M\otimes T^*M)$.

 \begin{proposition}\label{invcon}
$g\in H_2^p(M,T^*M\otimes T^*M)$ is a Riemannian metric on $M$ with $p>n/2$. 
If $\gt=\psi^{\frac{4}{n-2}}g$ is a conformal metric to  $g$ with $\psi\in H_2^p(M)$ and $\psi>0$ then $L$ is weakly conformally invariant, which means that
\begin{equation}\label{laplinvfai}
 \forall u\in H_1(M)\qquad \psi^{\frac{n+2}{n-2}}L_{\gt}(u)=L_g(\psi u)\quad weakly 
\end{equation}
Moreover if $\mu(g)>0$ then the conformal Laplacian $L_g=\Delta_g+\frac{n-2}{4(n-1)}R_g$ is invertible and coercive.
\end{proposition}

\begin{proof}
Recall that $\di v_{\gt}=\psi^{\Nn}\di v$ and
$$\forall u,w\in L^2(M)\quad (u,w)_{g,L^2}=\int_Muw\di v_g$$ is the scalar product in  $L^2(M)$ with the metric $g$.\\ 
For all $u,w\in H_1(M)$: 
\begin{equation*}
\begin{split}
(\psi^{\Nn}L_{\gt} u,w)_{g,L^2} &=(L_{\gt} u,w)_{\gt,L^2}\\
&=\int_M \gt(\nabla u ,\nabla w)+\an R_{\gt} uw \di v_{\gt}\\
                              &=\int_M \psi^2g(\nabla u ,\nabla w)+\an R_{\gt}\psi^{\frac{n+2}{n-2}} (uw\psi)\di v_{g}		      
\end{split}
\end{equation*}
We know that the scalar curvatures $R_g$ and $R_{\gt}$ are related by Yamabe  equation \eqref{yamabe}, which is equivalent to 

$$L_g\psi=\an R_{\gt}\psi^{\frac{n+2}{n-2}}\quad weakly $$
then
$$(L_g\psi,uw\psi)_{g,L^2}=\an (R_{\gt}\psi^{\frac{n+2}{n-2}},uw\psi)_{g,L^2}$$
Hence
\begin{equation}
\begin{split}
(\psi^{\Nn}L_{\gt} u,w)_{g,L^2} &=\int_M \psi^2g(\nabla u, \nabla w)+g(\nabla\psi ,\nabla (uw\psi))+\an R_g \psi(uw\psi)\di v_{g}\\
			      & =\int_M g(\nabla (\psi u), \nabla (w\psi))+\an R_g (\psi u)(w\psi)\di v_{g}\\
			      & =(\psi L_g(\psi u),w)_{g,L^2} 
\end{split}
\end{equation}
We used the fact that  $u\psi$ and $w\psi$ belong to $H_1(M)$, indeed we have the the following Sobolev embedding 
$$H_2^p(M)\subset C^{1-[n/p],\beta}(M),\; H^p_1(M)\subset L^{\frac{pn}{n-p}}(M)\text{ and } H_1(M)\subset L^{\frac{2n}{n-2}}(M)$$
 
 Let us prove that  $L_g$ is invertible and coercive. Let $\lambda$ be the smallest eigenvalue of  $L_g$ with positive eigenfunction $\varphi\in H_1(M)$, then 
$$\lambda\|\varphi\|_2^2=(L_g\varphi,\varphi)_{g,L^2}=I_g(\varphi)\|\varphi\|_N^2\geq \mu(g)\|\varphi\|_N^2>0$$
hence $\lambda>0$. We conclude the result, by applying proposition \ref{delta fu}.

\end{proof}

\section{Yamabe conformal invariant}\label{invariance de mu}

In the case of smooth metrics, $\mu(g)$ is conformally invariant, which means that if $g$ and $\gt$ are two smooth conformal metrics then $\mu(g)=\mu(\gt)$. The next proposition shows that we can extend this property to metrics in $H^p_2$. 
\begin{proposition}\label{invconforme}
Let  $M$ be a smooth compact manifold of dimesion $n\geq 3$. Let $g$ and $\gt=\psi^{\frac{4}{n-2}}g$ be two metrics in $H^p_2$, with $\psi\in H^p_2(M)$ positive. if $p>n/2$  then
$$\mu(g)=\mu(\gt)$$ 
\end{proposition}

\begin{proof}
Let $u\in H_1(M)$ be test function for the Yamabe functional   $I_g$. Notice that $E(u)=(L_g(u),u)_{g,L^2}$. then $$I_{\gt}(u)=(L_{\gt}(u),u)_{\gt,L^2}\|u\psi\|_N^{-2}$$ 
Using proposition \ref{invcon}, we deduce that
$$I_{\gt}(u)=(L_g(\psi u), \psi u)_{g,L^2}\|u\psi\|_N^{-2}$$ 
Finally
\begin{equation}\label{IgIg}
 I_{\gt}(u)=I_g(\psi u)
\end{equation}
Which implies that $\mu(g)=\mu(\gt)$. So this invariant depends only on the conformal class   $[g]$ and the manifold  $M$.
\end{proof}

\section{Green Function}\label{glc}

\begin{definition}Let $(M,g)$ be a compact Riemannian manifold and $P$ be a point in $M$. We call $G_P$ the Green function on $P$ of the linear operator $L$, if it satisfies
$$LG_P=\delta_{P}   (\Longleftrightarrow \forall f\in C^{\infty}(M)\quad \langle G_P,Lf\rangle =f(P))$$
\end{definition}
Proposition \ref{grgr} shows the existence of such function for the operator $L=\Delta+h$ with  a positive continuous function $h$. Unfortunately,  the method used to construct this function doesn't work when $h$ belongs to $L^p(M)$.  This case holds for the conformal Laplacian operator  $L_g$, because $R_g\in L^p(M)$. But, using proposition \ref{green conf}, we construct this function and we obtain corollary \ref{green coroll}. 

\begin{proposition}\label{grgr}
Let  $h$ be a positive continuous and $P\in M$. $g$ is a metric satisfying assumption $(H)$. There exists a unique Green function $G_{P}$ for the operator $L=\Delta_g+h$ which satisfies  $LG_P=\delta_{P}$ and
\begin{itemize} 
\item[$(i)$] $G_{P}$ is smooth in $B_{P_0}(\delta)-\{P\}$
\item[$(ii)$] $G_{P}\in C^2(M-\{P\})$ 
\item[$(iii)$] There exists $c>0$ such that for any $Q\in M-\{P\}$,  $|G_{P}(Q)|\leq c d(P,Q)^{2-n} $
\end{itemize}
\end{proposition}
\begin{proof}
 $G_{P}$ is unique because  $L$ is invertible. In fact, if $\lambda$ is an eigenvalue of $L$ and $\varphi$ is a positive eigenfunction  associated to $\lambda$ then 
$$\lambda\|\varphi\|^2_2=(L\varphi,\varphi)_{L^2}=E(\varphi)>0$$
Hence $\lambda>0$. To conclude, it is sufficient to apply proposition \ref{delta fu}. For the existence of such function, we follow Aubin's \cite{Aubin} construction for the Laplacian, in the case of smooth metrics. We choose a decreasing positive smooth radial function $f(r)$, equal to $1$ for $r<\delta/2$ and zero for $r\geq\delta(M)$ the injectivity radius of $M$. We define the following functions

\begin{gather*} 
H(P,Q)=\frac{f(r)}{(n-2)\omega_{n-1}}r^{2-n}\mbox{ with }r=d(P,Q)\\
\Gamma^1(P,Q)=-L_QH(P,Q) \\
\forall i\in\mathbb{N^*}\qquad\Gamma^{i+1}(P,Q)=\int_M\Gamma^{i}(P,S)\Gamma^1(S,Q)\di v(S)
\end{gather*}
Then
\begin{equation*}
|\Gamma^1(P,Q)| \leq c d(P,Q)^{2-n}
\end{equation*}
We show that
\begin{equation*}
\forall i\geq 1\qquad|\Gamma^i(P,Q)| \leq \begin{cases}& cd(P,Q)^{2i-n}  \hspace{2.75cm}\text{ if }2i<n\\
& c(1+\log d(P,Q))  \hspace{2cm}\text{ if } 2i=n\\
& c   \hspace{4.7cm}\text{ if } 2i>n
\end{cases}
\end{equation*}
In the last case $\Gamma^i$ is continuous.\\
More details are given in Aubin' book \cite{Aubin}. \\
The Green function of  $L$ is given by
\begin{equation}\label{greeen}
G_{P}(Q)=H(P,Q)+\sum_{i=1}^{k}\int_M\Gamma^{i}(P,S)H(S,Q)\di v(S)+F_{P}(Q)
\end{equation} 
where $F_{P}$ satisfies
$$LF_{P}=\Gamma^{k+1}(P,\cdot)$$ 
We choose $k=[n/2]$, $\Gamma^{k+1}(P,\cdot)$ is continuous. Regularity theorem \ref{reg} implies that $F_P$ is $C^2$. \\
$(i)$ $L_gG_{P}=0$ in $B_{P_0}(\delta)-\{P\}$ and the metric is smooth on $B_{P_0}(\delta)$, regularity theorem assure that $G_{P}$ is smooth on $B_{P_0}(\delta)-\{P\}$, with $P\in M$.\\
$(ii)$ We have also $LG_{P}=0$ in $M-\{P\}$. We conclude that $G_P$ is $C^2$ in $M-\{P\}$.\\
$(iii)$ In the expression \eqref{greeen}, we notice that the leading term, in the neighborhood of $P$, is $H(P,Q)$, then for all $P\neq Q$, $$|G_P(Q)|\leq cd(P,Q)^{2-n}$$
\end{proof}
 
\begin{proposition}\label{green conf}
Let $g$ be a metric in $H^p_2(M,T^*M\otimes T^*M)$, $\gt=\psi^{\frac{4}{n-2}}g$ is conformal to $g$ with $\psi\in H^p_2(M)$ positive and $p>n/2$. We suppose that $L_{\gt}$ admits a Green function on $P$, denoted $\tilde G_{P}$, then $L_g$ admits a Green function, denoted $G_P$ and it is given by
$$\forall Q\in M-\{P\}\qquad G_P(Q)=\psi(P)\psi(Q)\tilde G_P(Q)$$ 
\end{proposition}

\begin{proof}

For any function $\varphi\in C^{\infty}(M)$:
\begin{equation*}
\begin{split}
\langle \psi(P)\psi\tilde G_P,L_{g}\varphi\rangle _g &=\psi(P)\int_M\tilde G_P\psi L_{g}[\psi(\frac{\varphi}{\psi})]\di v_g\\
							  &=\psi(P)\int_M \tilde G_PL_{\gt}\frac{\varphi}{\psi}
								\di v_{\gt}\\
							  &=\psi(P)\langle\tilde G_P,L_{\gt}\frac{\varphi}{\psi}\rangle  _{\gt}\\
							  &=\varphi(P)
\end{split}
\end{equation*}
The second equality above is obtained by the weak conformal invariance of the conformal Laplacian (see proposition \ref{invcon}). We know that for any $Q\in M-\{P\}$
$$|\tilde G_P(Q)|\leq  cd(P,Q)^{2-n}$$
then $G_P\in L^s(M)$, for any $s\in [1,n/(n-2))$ and $L_{\gt}\frac{\varphi}{\psi}\in L^p(M)$ with $p>n/2$. We choose $s$ such that $\langle\tilde G_P,L_{\gt}\frac{\varphi}{\psi}\rangle  _{\gt}$ is finite. Hence the third equality is well defined.  
\end{proof}

\begin{corollary}\label{green coroll}
$g$  is a Riemannian metric, satisfying assumption $(H)$. If $\mu(g)>0$ then the conformal Laplacian $L_g$ admits a Green function $G_{P_0}$ which satisfies $LG_{P_0}=\delta_{P_0}$ and
\begin{itemize} 
\item[$(i)$] $G_{P_0}$ is smooth in  $B_{P_0}(\delta)-\{P_0\}$
\item[$(ii)$] $G_{P_0}\in H^p_2(M-B_{P_0}(r))$ for any $r>0$. 
\item[$(iii)$] There exists $c>0$ such that for any $Q\in B_{P_0}(\delta)-\{P_0\}$,  $|G_{P_0}(Q)|\leq c d(P_0,Q)^{2-n} $
\end{itemize}
\end{corollary}

\begin{proof}
$\mu(g)>0$,  $L_g$ is invertible. We deduce that $L_g$ admits a unique Green function. Using standard variational method (see the  proof of poroposition \ref{delta fu} ), we can show that the equation
\begin{equation}\label{eygiii}
\Delta_g \psi+ \frac{n-2}{4(n-1)}R_g\psi= \mu_{q,G}(g) \psi^{q-1}
\end{equation}
admits a positive solution $\psi\in H^p_2(M)$ when $2\leq q<N$, with 

$$ \mu_{q,G}(g)=\inf_{\psi\in H_1(M)-\{0\}} \frac{E(\psi)}{\|\psi\|_q^2}$$
Moreover, $g$ is smooth in $B_{P_0}(\delta)$, regularity theorem shows that  $\psi$ is also smooth in the same ball. The metric $\gt:=\psi^{\frac{4}{n-2}}g$ satisfies assumption $(H)$. Using  Yamabe equation \eqref{yamabe},  we deduce that the scalar curvature of  $\gt$  is $$R_{\gt}=\frac{4(n-1)}{n-2}\mu_{q,G}(g) \psi^{q-N}$$ 
Hence $R_{\gt}$ is positive continuous because $\mu_{q,G}(g)>0$. Now, we are able to use proposition \ref{grgr}, which assure the existence of the Green function  $\tilde G_{P_0}$ for $L_{\gt}$ with the metric $\gt$. using proposition \ref{green conf}, we conclude that $G_{P_0}=\psi(P_0)\psi\tilde G_{P_0}$ is the Green function of $L_g$. The metrics $g$ and $\gt$ are smooth in $B_{P_0}(\delta)$ and $\tilde G_{P_0}$ satisfies the properties of proposition \ref{grgr}, then the properties announced for $G_{P_0}$ are valid.
\end{proof}

\section{Existence theorem}

\begin{theorem}\label{conj aub}
Let $M$ be a smooth compact manifold of dimension $n\geq 3$, $g$ is a Riemannian metric which satisfies the assumption $(H)$. If $(M,g)$ is not conformal to the sphere $(S_n, g_{can})$ then $\mu(g)<K^{-2}(n,2)$.
\end{theorem}
This theorem assure that Aubin's conjecture \ref{Aubincon} still valid for any metric satisfying the assumption  $(H)$.  \\
To prove this theorem, we use the results of Aubin and Schoen, when the metric $g$ is smooth. The strategy is the following: we construct a test function for the functional $I_g$, with a support in small geodesic ball. Then the problem is local. We know that the metric $g$ is smooth in $B_{P_0}(\delta)$, so the proof of this theorem is the same as when the metric is smooth everywhere (this is the point where we need the assumption : $g$ is smooth in $B_{P_0}(\delta)$). After, we consider Aubin and Schoen's test functions.\\  
We need also the following result obtained by Aubin \cite{Aub2}, for the Green function of $L_g$: 
\begin{theorem}\label{aubm}
If $g$ is a  Cao--Günther metric, $L_g$ is invertible and the normalized Green function $G_{P_0}$ have the following expression
$$G_{P_0}(Q)=r^{2-n}+A+O(r)$$ 
in a neighborhood of $P_0$ with $r=d(P_0,Q)$, then $A>0$, except if $(M,g)$ is conformal to $(S_n,g_{can})$ for which $A=0$.
\end{theorem}

\begin{proof}[\textbf{Proof of theorem \ref{conj aub}}]
If $\mu(g)\leq 0$ then the inequality is obvious. From now until the end of the proof, we suppose that $\mu(g)>0$. without loss of generalities, we suppose that $g$ is a Cao--Günther metric given in theorem  \ref{caogun}. In fact, $\mu(g)$ is conformally invariant (see proposition \ref{invcon}).\\
There are two cases which can happen :\\
$(a)$ The case  $(M,g)$ is not conformally flat in a neighborhood of $P_0$ and $n\geq 6$. We define $\varphi_\e=\eta v_\e$, $\eta$ is a cut-off function with support in $B_{P_0}(2\e)$, $\eta=1$ in $B_{P_0}(\e)$, $2\e<\delta$ and
$$v_\e(Q)=\biggl(\frac{\e}{r^2+\e^2}\biggr)^{\frac{n-2}{2}}\quad r=d(P_0,Q)$$
$supp\varphi\subset B_{P_0}(\delta)$ and the metric $g$ is smooth in this ball, we obtain the following lemma (see Aubin \cite{Aub}): 
\begin{lemma}
\begin{equation*}
\mu(g)\leq I_g(\varphi_\e)\leq \begin{cases}& K^{-2}(n,2)-c|W_g(P_0)|^2\e^4+o(\e^4)\text{ si }n>6\\
& K^{-2}(n,2)-c|W_g(P_0)|^2\e^4\log \frac{1}{\e}+O(\e^4)\text{ si }n=6
\end{cases}
\end{equation*}
where $|W_g(P_0)|$ is the norm of the Weyl tensor on $P_0$.
\end{lemma}
 ( Lee et Parker \cite{LP} gave a simple proof of this lemma, using the conformal normal coordinates on $P_0$).
Using this lemma, we conclude that $\mu(g)<K^{-2}(n,2)$.\\
$(b)$ The case $(M,g)$ is conformally flat in a neighborhood of $P_0$ or $n=3,\; 4\text{ or }5$. In this coordinates system, the Taylor expansion of the Green function is:
$$G_{P_0}(Q)=r^{2-n}+A+O(r)$$
with $r=d(P_0,Q)$ (see Lee and Parker's paper \cite{LP} for the proof of this expansion).\\
If $g$ satisfies assumption $(H)$ and $(M,g)$ is not conformal to $(S_n,g_{can})$, then theorem \ref{aubm} assure that $A>0$. 
Hence we can consider Shoen's test function  $\varphi_\e$, defined for any  $Q\in M$ by:
\begin{equation*}
\varphi_\e(Q)= \begin{cases}& v_\e(Q)\text{ if }Q\in B_{P_0}(\rho_0)\\
& \e_0[G_{P_0}-\eta (G_{P_0}-r^{2-n}-A)](Q) \text{ if }Q\in B_{P_0}(2\rho_0)-B_{P_0}(\rho_0)\\
& \e_0 G_{P_0}(Q)\text{ if }Q\in M-B_{P_0}(2\rho_0)
\end{cases}
\end{equation*}
with $2\rho_0<\delta$, $(\frac{\e}{\rho_0^2+\e^2})^{(n-2)/2}=\e_0(\rho_0^{2-n}+A)$ and 
$\eta$ is a smooth nonnegative decreasing function on $\mathbb{R}_+$, with support in $(-2\rho_0,2\rho_0)$, equal to $1$ in $[0,\rho_0]$, the gradient 
$|\nabla\eta(r)|\leq\rho_0^{-1}$. $g$ is smooth in $B_{P_0}(2\rho_0)\subset B_{P_0}(\delta)$ and  $G_{P_0}\in H^p_2(M-B_{P_0}(\rho_0))$ (see corollary \ref{green coroll}), then we have the estimate of  $\mu(g)$, obtained by  Schoen\cite{Schoen}:
\begin{lemma}
 $$\mu(g)\leq I_g(\varphi_\e)\leq K^{-2}(n,2)+c\e_0^2(c\rho_0-A)$$
\end{lemma}
The fact that $A>0$ allows us to choose  $\rho_0$ sufficiently  small ($c\rho_0<A$) such that $\mu(g)<K^{-2}(n,2)$.\\
\end{proof}

Now, we can state the main theorem which solves the problem \ref{yam sing} for any metric which satisfies assumption $(H)$. 

\begin{theorem}\label{inegg}
Let  $M$ be a smooth compact manifold of dimension $n\geq 3$ and $g$ be a metric satisfying assumption  $(H)$.  
There exists a metric $\tilde g$ conformal to $g$ such that the scalar curvature $R_{\tilde g}$ is constant everywhere. This metric solves the problem \ref{yam sing}.
\end{theorem}

It means that we can always solve the equation of type Yamabe \eqref{AF} when $h=\frac{n-2}{4(n-1)}R_g$.

\begin{proof}

If $(M,g)$ is conformal to $(S_n, g_{can})$ then the result is obvious because the scalar curvature of $(S_n, g_{can})$ is constant. Otherwise $(M_n,g)$ is not conformal to $(S_n,g_{can})$. In this case, we have the inequality 
$$\mu(g)<K^{-2}(n,2)$$ 
given by theorem \ref{conj aub}. Using theorem \ref{cg}, we get a positive solution $\psi\in H^p_2(M)$ of \eqref{AF}, where $h=\frac{n-2}{4(n-1)}R_g$ and $\tilde h=\mu(g)$. Using Yamabe equation \eqref{yamabe}, we deduce that the metric $\gt=\psi^{\frac{4}{n-2}}g$ has a constant scalar curvature $R_{\gt}=\frac{4(n-1)}{n-2}\mu(g)$.
\end{proof}

\section{Uniqueness of solutions}

When the metrics are smooth, if  $\mu(g)$ is nonpositive then the solutions of the Yamabe equation \eqref{yamabe} are proportional. The following theorem generalizes the uniqueness theorem in the singular case.  

\begin{theorem}\label{unique}
 Let $g$ be a metric in   $H^p_2(M,T^*M\otimes T^*M)$, with $p>n$. If $\mu(g)\leq 0$ then the solutions of \eqref{yamabe} are proportional.
\end{theorem}

\begin{proof}
Let $\varphi_1$ and $\varphi_2$ two positive solutions of \eqref{yamabe}. The metrics $g_i=\varphi_i^{\frac{4}{n-2}}g$ have a constant scalar curvatures $R_i$, where $i=1$ or $2$. Define $\psi=\frac{\varphi_1}{\varphi_2}$, then $g_1=\psi^{\frac{4}{n-2}}g_2$. It implies that   $\psi$ satisfies
\begin{equation}\label{unieq}
\Delta_{g_2} \psi+ \frac{n-2}{4(n-1)}R_{2}\psi=  \frac{n-2}{4(n-1)}R_{1} \psi^{\frac{n+2}{n-2}} 
\end{equation}

By  regularity theorem \ref{reg}, we deduce that  $\psi$ is $C^{2,\beta}$ because the coefficient of the Laplacian are $C^0$. In fact, in a local coordinates system :
$$\Delta_g\psi=-\nabla_i\nabla^i\psi=-g^{ij}(\partial_{ij}\psi-\Gamma^k_{ij}\partial_k\psi)$$
and the Christoffels are in $H^p_1(M)$ then continuous if $p>n$. In other hands, notice that    $R_1$, $R_2$ have the same sign.  Hence, if $\mu(g)<0$ then $R_i<0$ for $i=1$ and 2. Let  $Q_1\in M$ (resp.  $Q_2\in M$) be a point for  which $\psi$ is maximal  (resp. minimal ).  Then $\Delta_{g_2}\psi(Q_1)\geq 0$ and  $\Delta_{g_2}\psi(Q_2)\leq 0$. Hence, if we evaluate  equation \eqref{unieq} at  $Q_1$ and $Q_2$, we obtain :
$$\psi^{\frac{4}{n-2}}(Q_1)\leq \frac{R_2}{R_1}\mbox{ and }\psi^{\frac{4}{n-2}}(Q_2)\geq \frac{R_2}{R_1}$$ 
We conclude that $\psi=\frac{R_2}{R_1}$, $\varphi_1$ and $\varphi_2$ are proportional.\\
If $\mu(g)=0$ then $R_1=R_2=0$ and \eqref{unieq} becomes $\Delta_{g_2}\psi=0$, hence $\psi$ is constant.
\end{proof}



\end{document}